\documentclass[12pt]{amsproc} 
\usepackage{euscript}
\newcommand\Diff{\operatorname{Diff}}
\newcommand{\KKK}{\EuScript{K}}
\newcommand{\QQ}{\mathbb{Q}}

\newcommand{\ZZ}{\mathbb{Z}}
\newcommand{\PP}{\mathbb{P}}
\newcommand{\NN}{\mathbb{N}}

\newcommand{\Supp}{\operatorname{Supp}}
\newcommand{\Center}{\operatorname{Center}}
\newcommand{\down}[1]{\left\lfloor #1\right\rfloor}

\newtheorem*{theorem}{Theorem}
\newtheorem*{proposition}{Proposition}
\newtheorem*{lemma}{Lemma}
\theoremstyle{definition}
\newtheorem*{example}{Example}

\date{}
\title{A remark on resolution of terminal singularities}

\author{Yu.~G.~Prokhorov}

\thanks{The author was
partially supported by the INTAS-OPEN  269 and Russian Foundation
of Basic Research  99-01-01132}

\address{
Department of Algebra, Faculty of Mathematics, Moscow State
Lomonosov University, Moscow 117234, Russia}
\email{prokhoro@mech.math.msu.su}

\begin{document}
\maketitle
\begin{abstract}
Let $(Z,o)$ be a three-dimensional terminal singularity of type
$cA/r$. We prove that all exceptional divisors over $o$ with
discrepancies $\le 1$ are rational.
\end{abstract}

Let $(Z,o)$ be a three-dimensional terminal singularity of index
$r\ge 1$ and let $\varphi\colon \tilde Z\to Z$ be a resolution.
Let $S\subset \tilde Z$ be an exceptional divisor such that
$\Center(S)=o$. It is clear that the birational type of $S$ does
not depend on $g$. By \cite[2.14]{Reid-canonical} the surface $S$
is birationally ruled. We say that the corresponding discrete
valuations $\nu=\nu_S$ of the function field $\KKK(Z)$ is
\emph{rational} (resp. \emph{birationally ruled}) if so is the
surface $S$. We are interested in the existence of rational
exceptional divisors over $o\in Z$ with small discrepancies:

\begin{theorem}
Let $(Z, o)$ be a three-dimensional terminal singularity of type
$cA/r$, $r\ge 1$ and let $\nu$ be a divisorial discrete valuation
of the function field $\KKK(Z)$ such that $a(\nu)\le 1$ and
$\Center_Z(\nu)=o$. Then $\nu$ is rational.
\end{theorem}

Recall that according to the classification \cite{Mori},
\cite{RY}, $(Z,o)$ belongs to one of the following classes:
$cA/r$, $cAx/4$, $cAx/2$, $cD/2$, $cD/3$, $cE/2$, $cD$, $cE$. It
was proved in the series of works \cite{Ka}, \cite{Mar},
\cite{Sh-s} (see also \cite{Sh2}, and \cite{Ha}) that for any
$i\in \NN\setminus r\NN$ there exists an exceptional divisor $S$
with center at $o$ and discrepancy $a(S)=i/r$. Exceptional
divisors with discrepancies $<1$ appear in any resolution.
Divisors with discrepancy $1$ and $\Center=o$ appear in any
\emph{divisorial} resolution, i.e., in a resolution such that the
exceptional set is of pure codimension $1$.

\begin{proof}[Proof of Theorem]
Let $F\in |{-}K_Z|$ be a general member. Then $(F, o)$ is a Du Val
singularity (of type $A_n$) \cite[6.4]{RY}. By the Inversion of
Adjunction \cite[17.6]{Ut} the pair $(Z,F)$ is plt. Let $q\colon Z
^{q} \to Z $ be $\QQ$-factorialization (see, e.g.,
\cite[6.7]{Ut}). Then $Z^{q}$ has (terminal) singularities of
types $cA/r_i $. Indeed, the surface $F^q: = q^{-1} (F)$ contains
all singular points of $Z^{q} $ and $F^q\in | -K_{Z^{q}} | $.
Since $F^q\to F $ is a crepant morphism, $F^q $ has only
singularities of types $A_{n_i}$. Thus, replacing $Z$ with
$Z^{q}$, we may assume that $Z$ is $\QQ$-factorial.

Let $S$ be an exceptional divisor with center at $o$ and
discrepancy $a(S)\le 1$. Then $a(S,F)<1$. Since $K_Z+F$ is
linearly trivial, we have $a(S,F)\le 0$. Further, there is an
$1$-complement $K_F+\Theta$ of $K_F$ near $o$ (see
\cite[5.2.3]{Sh}). According to \cite[4.4.1]{Lect} this complement
can be extended to $Z$, i.e., there is an (integral) Cartier
divisor $F'$ such that $F'|_F=\Theta$, $K_Z+F+F'\sim 0$, and
$(Z,F+F')$ is lc. Then $a(S,F+F')=-1$. Now our theorem is a
consequence of the following simple fact.
\end{proof}

\begin{proposition}
Let $(Z, o)$ be a three-dimensional $\QQ$-factorial terminal
singularity and let $\nu $ be a divisorial discrete valuation of
the field $\KKK (Z)$. Assume that there is a boundary $D $ such
that the pair $(Z, D)$ is lc and $a(\nu, D) = -1$. Then
\par\smallskip\noindent
{\rm (i)} The valuation $\nu $ is rational or is
birationally a ruled surfaces over an elliptic curve;
\par\smallskip\noindent
{\rm (ii)} if, moreover, $\down D $ there are at least two
components passing through $o $, then $\nu $  is rational.
\end{proposition}

\begin{proof}
According to \cite[17.10] {Ut} there is a blowup $f\colon X\to Z $
with irreducible exceptional divisor $S$ representing the
valuation $\nu $ such that the log divisor $K_X+S+D_X=f^* (K_Z+D)
$ is lc. Here $D_X $ is the proper transform of $D $. In this
situation we have $\rho (X/Z) =1 $. Therefore, $D_X\equiv -
(K_X+S)$ is $f$-ample.

Consider a minimal log terminal modification $g\colon Y\to X$ of
the pair $(X,S+D_X)$ (see, e.g., \cite[3.1.3]{Lect}), i.e., a
blowup such that $Y$ is $\QQ$-factorial and
\[
K_Y+S_Y+D_Y+E=g^*(K_X+S+D_X),
\]
is dlt. Here $S_Y$, $D_Y$ are proper transforms of $S$, $D_X$,
respectively, and $E=\sum E_i$ is a (reduced) exceptional divisor
with $a(E_i,S+D_X)=-1$. Denote $\Delta:=\Diff_{S_Y}(D_Y+\sum E_i)$
and $\Omega:=f^*D_X|_{S_Y}$. By \cite[17.7]{Ut} the surface $S_Y$
is normal and the pair $(S_Y,\Delta)$ is lc. The assertion of (i)
easily follows by the lemma below.

To prove (ii) we assume that $S_Y$ is not rational and $\down D$
is not irreducible. Then $S\cap \down{D_X}$ has at least two
irreducible components. So is $S_Y\cap \down{D_Y+E}$. By the lemma
below $\down{\Delta}$ has exactly two components (contained in
$\down{D_Y+E}$) and the pair $(S_Y,\Delta)$ is plt. Further,
$(S_Y,\Delta-\varepsilon\Omega)$ is klt whenever $0<\varepsilon$.
For $0<\varepsilon\ll 1$, the pair
$(S_Y,\Delta-\varepsilon\Omega)$ is a klt log del Pezzo (because
$\Omega$ is nef and big). In this situation, $S_Y$ is rational
(see, e.g., \cite[5.4.1]{Lect}).
\end{proof}

\begin{lemma}
Let $(S, \Delta)$ be a projective log surface with $\kappa (S) =
-\infty $. Assume that $K_X +\Delta $ is lc and numerically
trivial and the surface $S $ is nonrational. Then  $S $  is
birationally ruled over an elliptic curve and there exists at most
two divisors with discrepancy $a (\ , \Delta) = -1$.
\end{lemma}
\begin{proof} [Proof (cf. {\cite[6.9] {Sh}})]
Replace $S$ with its minimal resolution and $\Delta $ with its
crepant pull-back. There is a contraction $\phi\colon S\to C $
(with general fiber $\PP^1$) onto a curve $C$ of genus $p_a (C)
\ge 1 $. In this situation pair $(S, \Delta)$ has only canonical
singularities and all components of $\Delta$ are horizontal
\cite[8.2.2-8.2.3] {Lect}. Hence, divisors with discrepancy $a (\
, \Delta) = -1 $ are exactly components of $\down {\Delta} $. As
an immediate consequence we have that the number of divisors with
discrepancy $-1 $ is at most two. If $\down {\Delta} \neq 0 $,
then for any component $\Delta_i\subset \down {\Delta} $ we have
$2p_a (\Delta_i) -2\le (K_S +\Delta) \cdot \Delta_i=0 $.
Therefore, $p_a (C) \le p_a (\Delta_i) \le 1 $. It remains to
consider the case, when the pair $(S, \Delta)$ is klt. Again for
any component $\Delta_i\subset \Supp (\Delta)$ we have $
\Delta_i^2\le 0$ (otherwise $(S,\Delta+\varepsilon\Delta_i)$ is a
klt log del Pezzo). As above, $p_a (C) \le p_a (\Delta_i) =
\frac12 ( K_S +\Delta_i) \cdot \Delta_i+1\le \frac12 (K_S +\Delta)
\cdot \Delta_i+1=1$.
\end{proof}

Note that the assertion (ii) of our theorem is not true for other
types of terminal singularities:

\begin{example}[\cite{Ka}]
Let $(Z,o)$ be a terminal $cAx/2$-singularity
$\{x^2+y^2+z^{4m}+t^{4m}=0\}/\ZZ_2(0,1,1,1)$. Consider the
weighted blowup with weight $\frac12(2m,2m+1,1,1)$. Then the
exceptional divisor $S$ is given in $\PP(2m,2m+1,1,1)$ by the
equation $x^2+z^{4m}+t^{4m}=0$. It is reduced, irreducible, and
$a(S)=1/2$. It is easy to see that $S$ is a birationally ruled
surface over a hyperelliptic curve of genus $2m-1$.
\end{example}

\end{document}